\documentclass[12pt,a4paper,reqno]{amsart}

\usepackage{amssymb,mathtools}

\edef\restoreparindent{\parindent=\the\parindent\relax}
\usepackage{parskip}
\restoreparindent

\makeatletter
\@namedef{subjclassname@2020}{\textup{2020} Mathematics Subject Classification}
\makeatother

\newtheorem{thm}{Theorem}
\newtheorem{lem}{Lemma}

\newtheorem{cor}{Corollary}
\newtheorem{conj}{Conjecture}

\newtheorem{Thm}{Theorem}

\newtheorem{Lem}{Lemma}

\theoremstyle{definition}
\newtheorem{rem}{Remark}
\newtheorem{defn}{Definition}
\newtheorem{exm}{Example}

\newcommand{\IR}{{\mathbb R}}

\newcommand{\IC}{{\mathbb C}}
\newcommand{\ID}{{\mathbb D}}

\newcommand{\IT}{{\mathbb T}}

\newcommand{\real}{{\operatorname{Re}\,}}
\newcommand{\ima}{{\operatorname{Im}\,}}
\newcommand{\eit}{{e^{i\theta}}}

\newcommand{\be}{\begin{equation}}
\newcommand{\ee}{\end{equation}}

\newcommand{\blem}{\begin{lem}}
\newcommand{\elem}{\end{lem}}

\newcommand{\bdefn}{\begin{defn}}
\newcommand{\edefn}{\end{defn}}

\newcommand{\bthm}{\begin{thm}}
\newcommand{\ethm}{\end{thm}}

\newcommand{\bcor}{\begin{cor}}
\newcommand{\ecor}{\end{cor}}

\newcommand{\bconj}{\begin{conj}}
\newcommand{\econj}{\end{conj}}

\newcommand{\brem}{\begin{rem}}
\newcommand{\erem}{\end{rem}}

\newcommand{\bpf}{\begin{proof}}
\newcommand{\epf}{\end{proof}}

\begin{document}
	
\bibliographystyle{abbrv}

\title[Note on real and imaginary parts of harmonic quasiregular mappings]
{Note on real and imaginary parts of harmonic quasiregular mappings}

\author{Suman Das}
\address{Suman Das\vskip0.05cm Department of Mathematics with Computer Science, Guangdong Technion - Israel
	Institute of Technology, Shantou, Guangdong 515063, P. R. China.}
\email{suman.das@gtiit.edu.cn}

\author{Antti Rasila}
\address{Antti Rasila \vskip0.05cm Department of Mathematics with Computer Science, Guangdong Technion - Israel
	Institute of Technology, Shantou, Guangdong 515063, P. R. China. \vskip0.025cm Department of Mathematics, Technion - Israel
	Institute of Technology, Haifa 3200003, Israel.}
\email{antti.rasila@gtiit.edu.cn; antti.rasila@iki.fi}



\begin{abstract}
If $f=u+iv$ is analytic in the unit disk $\ID$, it is known that the integral means $M_p(r,u)$ and $M_p(r,v)$ have the same order of growth. This is false if $f$ is a (complex-valued) harmonic function. However, we prove that the same principle holds if we assume, in addition, that $f$ is $K$-quasiregular in $\ID$. The case $0<p<1$ is particularly interesting, and is an extension of the recent Riesz type theorems for harmonic quasiregular mappings by several authors. Further, we proceed to show that the real and imaginary parts of a harmonic quasiregular mapping have the same degree of smoothness on the boundary.
\end{abstract}

\maketitle
\pagestyle{myheadings}
\markboth{S. Das and A. Rasila}{Note on real and imaginary parts of harmonic quasiregular mappings}



\section{Introduction and background}\label{sec1}

\subsection{Notations and preliminaries}
Let $\ID$ denote the open unit disk in the complex plane and $\IT$ be the unit circle. For a function $f$ analytic in $\ID$,
the \textit{integral means} are defined as $$M_p(r, f) \coloneqq \left( \frac{1}{2\pi}\int_{0}^{2\pi}\vert f(r e^{i\theta}) \vert ^p\, d\theta \right)^{1/p} \quad \text{if} \quad 0 < p < \infty,$$ and $$M_\infty(r,f) \coloneqq \sup_{ \vert z \vert=r} \vert f(z) \vert.$$
It is said that $f$ is in the \textit{Hardy space} $H^p$ $(0 < p \le \infty)$ if $$\|f\|_p \coloneqq \lim_{r \to 1^-} M_p(r,f)<\infty.$$ A function $f \in H^p$ has the radial limit $$f(e^{i\theta})\coloneqq \lim_{r \to 1^-} f(r\eit)$$ in almost every direction, and $f(\eit)\in L^p(\IT)$.
Detailed surveys on Hardy spaces and integral means can be found in, for example, the book of Duren \cite{Duren}. Throughout this paper, we follow notations from \cite{Duren}.

A complex-valued function $f=u+iv$ is harmonic in $\ID$ if $u$ and $v$ are real-valued harmonic functions in $\ID$. Every such function has a unique representation $f=h+\overline{g}$, where $h$ and $g$ are analytic in $\mathbb{D}$ with $g(0)=0$. Analogous to the $H^p$ spaces, the \textit{harmonic Hardy spaces} $h^p$ are the class of harmonic functions $f$ satisfying $\|f\|_p < \infty.$ 

\subsection{Growth of conjugate functions}
Given a real-valued harmonic function $u$ in $\ID$, let $v$ be its harmonic conjugate with $v(0)=0$. It is a natural question that if $u$ has a certain property, whether so does $v$. In the context of boundary behaviour, this is answered by a celebrated theorem of M.~Riesz.
\begin{Thm}\label{MRiesz}\cite[Theorem 4.1]{Duren}
	If $u \in h^p$ for some $p$, $1<p<\infty$, then its harmonic conjugate $v$ is also of class $h^p$. Furthermore, there is a constant $A_p$, depending only on $p$, such that $$M_p(r,v) \le A_p \, M_p(r,u),$$ for all $u \in h^p$.
\end{Thm}
Curiously, the theorem fails for $p=1$ and $p=\infty$, examples can be found in \cite[p. 56]{Duren}. Although the harmonic conjugate of an $h^1$-function need not be in $h^1$, Kolmogorov proved that it does belong to $h^p$ for all $p<1$. Later, Zygmund established that the condition $|u| \log^+ |u| \in L^1(\IT)$ is the ``minimal" growth restriction on $u$ which implies $v\in h^1$. We refer to the paper of Pichorides \cite{Pichorides} for the optimal constants in the Riesz, Kolmogorov, and Zygmund theorems.

In  \cite{HL2}, Hardy and Littlewood showed that in the case $0<p<1$, Riesz's theorem is false in a much more comprehensive sense. Kolmogorov's result might suggest that if $u\in h^p$, then $v$, while not necessaily in $h^p$, should belong to $h^q$ for $0<q<p$. But this is false, and in fact, $v$ need not belong to $h^q$ for any $q>0$.

Nevertheless, They proved that the symmetry is restored in these latter cases if instead of the boundedness of the means, one considers their order of growth.

\begin{Thm}\label{conj}\cite[Theorem 4]{HL2}
	Let $0<p\le \infty$ and $\beta>0$. Suppose $f=u+iv$ is analytic in $\ID$, and $$M_p(r,u) = O \left(\frac{1}{(1-r)^{\beta}}\right).$$ Then $$M_p(r,v) = O \left(\frac{1}{(1-r)^{\beta}}\right).$$
\end{Thm}

The proof of this theorem is based on an extremely complicated (as remarked by the authors themselves) result, which can be stated as follows.

\begin{Thm}\label{hl_new}\cite[Theorems 2 \& 3]{HL2}
	If $f=u+iv$ is analytic in $\ID$, and $$M_p(r,u) = O \left(\frac{1}{(1-r)^{\beta}}\right),\quad 0<p\le \infty,\quad \beta\ge 0,$$ then $$M_p(r,f') = O \left(\frac{1}{(1-r)^{\beta+1}}\right).$$ Further, the converse is true for all $\beta>0$.
\end{Thm}
Let us note that the functions $|u|^p$ and $|v|^p$ are subharmonic when $p\ge 1$, but not when $p<1$, and therefore, $M_p(r,u)$ and $M_p(r,v)$ are not necessarily monotonic for $p<1$. This is the principal difficulty in dealing with the case $0<p<1$ for harmonic functions.


\subsection{Riesz theorem for harmonic quasiregular mappings}
For $K\ge 1$, a sense-preserving harmonic function $f=h+\overline{g}$ is said to be $K$-\textit{quasiregular} if its \textit{complex dilatation} $\omega = g'/h'$ satisfies the inequality $$|\omega(z)| \le k <1 \quad (z\in \ID),$$ where \be\label{eq1}k\coloneqq\frac{K-1}{K+1}.\ee The function $f$ is $K$-\textit{quasiconformal} if it is $K$-quasiregular as well as homeomorphic in $\ID$. One can find the $H^p$-theory for quasiconformal mappings in, for example, the paper of Astala and Koskela \cite{ast_kos}. It is worth mentioning that harmonic quasiconformal mappings have generated considerable interest in recent times, perhaps from a novel point of view. In \cite{wang_rasila}, Wang et al. constructed independent extremal functions for harmonic quasiconformal mappings, which were then further explored by Li and Ponnusamy in \cite{Li_Pon}. Later in \cite{DHR2}, it was established that the newly constructed functions are also extremal for Baernstein type results in the Hardy space of harmonic quasiconformal mappings.

Suppose $f=u+iv$ is a harmonic function in $\ID$, and $u\in h^p$ for some $p>1$. Then the imaginary part $v$ does not necessarily belong to $h^p$, i.e., Riesz theorem is not true for harmonic functions. One naturally asks under which additional condition(s) a harmonic analogue of the Riesz theorem would hold. Recently, Liu and Zhu \cite{Liu_Zhu} showed that such a condition is the quasiregularity of $f$.
\begin{Thm}\label{LZ}\cite{Liu_Zhu}
	Let $f=u+iv$ be a harmonic $K$-quasiregular mapping in $\ID$ such that $u \ge 0$ and $v(0)=0$. If $u \in h^p$ for some $p\in (1,2]$, then also $v$ is in $h^p$. Furthermore, there is a constant $C(K,p)$, depending only on $K$ and $p$, such that $$M_p(r,v) \le C(K,p)M_p(r,u).$$ Moreover, if $K=1$, i.e., $f$ is analytic, then $C(1,p)$ coincides with the optimal constant in the Riesz theorem.
\end{Thm}

The condition $u\ge 0$ was subsequently removed by Chen and Huang \cite{chen_huang}, who remarkably extended the result for all $p\in (1,\infty)$. Later in \cite{Kalaj_Kolmogorov}, Kalaj produced a couple of Kolmogorov type theorems for harmonic quasiregular mappings. Very recently, a quasiregular analogue of Zygmund's theorem has been obtained by Kalaj \cite{Kal_Zyg}, and also independently by the present authors and Huang \cite{DHR1}.

The purpose of this paper is to show that the real and imaginary parts of a harmonic quasiregular mapping have the same order of growth for all $p>0$. This extends Theorem \ref{LZ} to the cases $0<p<1$ and $p=\infty$. The main results and their proofs are presented in the next section.

\section{Main results and proofs}\label{sec_results}
In what follows, we always assume that $K$ and $k$ are related by \eqref{eq1}.
\bthm\label{thm1}
Suppose $0<p\le \infty$ and $\beta>0$, and let $f=u+iv$ be a harmonic $K$-quasiregular mapping in $\ID$. If $$M_p(r,u) = O \left(\frac{1}{(1-r)^{\beta}}\right),$$ then $$M_p(r,v) = O \left(\frac{1}{(1-r)^{\beta}}\right).$$
\ethm

\bpf
For $1<p< \infty$, we could apply the result of Chen and Huang \cite{chen_huang}, but here we shall give a simple proof which makes no appeal to this deeper result.

Let us write $f=h+\overline{g}$, and let $F=h+g$. Then $$\real F = \real f = u.$$ If $M_p(r,u)$ has the given order of growth,
it follows from Theorem \ref{hl_new} that $$M_p(r,F') = O \left(\frac{1}{(1-r)^{\beta+1}}\right).$$ Now, we observe $$F'=h'+g'=(1+\omega)h',$$ so that $$|F'|\ge (1-|\omega|)|h'| \ge (1-k)|h'|,$$ as $|\omega| \le k$. This readily implies $$M_p(r,h')\le \frac{1}{1-k}M_p(r,F') = O \left(\frac{1}{(1-r)^{\beta+1}}\right).$$ Since $|g'|\le k|h'|$, we also have $$M_p(r,g') = O \left(\frac{1}{(1-r)^{\beta+1}}\right).$$ Therefore, the converse part of Theorem \ref{hl_new} shows that $$M_p(r,h) = O \left(\frac{1}{(1-r)^{\beta}}\right)=M_p(r,g).$$ For $1 \le p \le \infty$, Minkowski's inequality gives $$M_p(r,f) \le M_p(r,h)+M_p(r, g),$$ while for $0<p<1$, we have $$M_p^p(r,f) \le M_p^p(r,h)+M_p^p(r, g).$$ In either case, we find that $$M_p(r,f) = O \left(\frac{1}{(1-r)^{\beta}}\right),$$ which, in turn, implies $$M_p(r,v) = O \left(\frac{1}{(1-r)^{\beta}}\right).$$ This completes the proof.
\epf

The next theorem deals with the case $\beta=0$. If $f=u+iv$ is harmonic $K$-quasiregular and $u\in h^p$ for some $p<1$, then of course, $v$ need not be in any $h^q$, as discussed before. Nevertheless, it is still possible to give an estimate on $M_p(r,v)$, as we show in Theorem \ref{thm_new}. The proof is somewhat similar to that of Theorem \ref{thm1}, and relies on the following lemma from \cite{DK3}.
\begin{Lem} \label{lem1}\cite{DK3}
	Let $0 < p < 1$. Suppose $f=h+\overline{g}$ is a locally univalent, sense-preserving harmonic function in $\ID$ with $f(0)=0$. Then
	$$
	\| f\|_p^p \le C \int_{0}^{1} (1-r)^{p-1} M_p^p(r,h')\, dr,
	$$
	where $C>0$ is a constant independent of $f$.
\end{Lem}

\bthm\label{thm_new}
Suppose $f=u+iv$ is a harmonic $K$-quasiregular mapping in $\ID$, and $u \in h^p$ for some $p \in (0,1)$. Then $$M_p(r,v) = O\left(\left(\log \frac{1}{1-r}\right)^{1/p}\right).$$
\ethm

\bpf
As before, we write $f=h+\overline{g}$ and $F=h+g$. Since $M_p(r,u)$ is bounded, an appeal to Theorem \ref{hl_new}, for $\beta=0$, shows that $$M_p(r,F') = O \left(\frac{1}{1-r}\right).$$ The quasiregularity of $f$, like in the previous proof, then implies $$M_p(r,h') = O \left(\frac{1}{1-r}\right).$$ Without any loss of generality, we assume that $f(0)=0$. For $0<r<1$, let $f_r(z)=f(rz)$. Applying Lemma \ref{lem1} for the function $f_r$, we find
\begin{align*}
M_p^p(r,f) & \le C\int_{0}^1 (1-t)^{p-1}M_p^p(rt,h')\, dt \le C\int_{0}^1 \frac{(1-t)^{p-1}}{(1-rt)^p}\, dt\\ & = C\left[\int_{0}^r \frac{(1-t)^{p-1}}{(1-rt)^p}\, dt+\int_{r}^1 \frac{(1-t)^{p-1}}{(1-rt)^p}\, dt\right]\\ & \le C \left[\int_{0}^r \frac{1}{1-t}\, dt+\frac{1}{(1-r)^p}\int_{r}^1 (1-t)^{p-1}\, dt \right]\\&= O\left(\log \frac{1}{1-r}\right).
\end{align*} Therefore, it follows that $$M_p(r,v) \le M_p(r,f) = O\left(\left(\log \frac{1}{1-r}\right)^{1/p}\right).$$ The proof is thus complete.
\epf

Generally speaking, Theorem \ref{thm1} suggests that the real and imaginary parts of a harmonic quasiregular mapping have the same ``order of infinity". We now wish to show that they also have the same degree of smoothness on the boundary (see Theorem \ref{thm2}).

Let $\Lambda_\alpha$ $(\alpha > 0)$ be the class of functions $\varphi: \IR \to \IC$ satisfying a H\"older condition of order $\alpha$, i.e., $$\vert\varphi(x)-\varphi(y)\vert \le A \vert x-y\vert^\alpha,$$ for some constant $A>0$. If $\alpha > 1$, $\Lambda_\alpha$ is the class of constant functions, hence we restrict attention to the case $0<\alpha\le 1$. Clearly, $\Lambda_{\beta}\subset\Lambda_{\alpha}$ for $\alpha < \beta$.

The following principle of Hardy and Littlewood says that an analytic function $f$ is H\"older continuous on the boundary if $f'$ has a `slow' rate of growth, and conversely.

\begin{Thm}\label{smooth}\cite[Theorem 40]{HL1}
	Let $f$ be an analytic function in $\ID$. Then $f$ is continuous in the closed disk $\overline{\ID}$ and $f(\eit) \in \Lambda_\alpha$ $(0<\alpha \le 1)$, if and only if $$\vert f'(r\eit)\vert = O\left(\frac{1}{(1-r)^{1-\alpha}}\right).$$
\end{Thm}

We are now prepared to discuss the final result of this paper.

\bthm\label{thm2}
Let $f=u+iv$ be a harmonic $K$-quasiregular mapping in $\ID$, and suppose $u$ is continuous in $\overline{\ID}$. If $u(\eit) \in \Lambda_{\alpha}$, $0<\alpha < 1$, then $v$ is continuous in $\overline{\ID}$ and $v(\eit) \in \Lambda_{\alpha}$.
\ethm

\bpf
First, we note that if $v$ is continuous on $\IT$, then $v(r\eit)$ is the Poisson integral of $v(\eit)$. Hence, the continuity of $v(\eit)$ would imply the continuity of $v$ in $\overline{\ID}$. Therefore, it is enough to show that $v(\eit) \in \Lambda_{\alpha}$.

Now, suppose $u(\eit)\in \Lambda_{\alpha}$ and $f=h+\overline{g}$. As before, we write $F=h+g$ so that $$\real F = \real f=u.$$ Since $u$ is continuous in $\overline{\ID}$, we can represent $F$ by the Poisson integral formula $$F(z) = \frac{1}{2\pi}\int_{0}^{2\pi} \frac{e^{it}+z}{e^{it}-z}\, u(e^{it})\, dt+i\ima F(0).$$ This implies \begin{align*}
	F'(z) &= \frac{1}{2\pi}\int_{0}^{2\pi}\frac{\partial}{\partial z} \left(\frac{e^{it}+z}{e^{it}-z}\right)\, u(e^{it})\, dt\\ & = \frac{1}{\pi}\int_{0}^{2\pi} \frac{e^{it}}{(e^{it}-z)^2}\, u(e^{it})\, dt.
\end{align*}
Therefore, for $z=r\eit$, we have \be\label{eq01}
F'(r\eit) = \frac{1}{\pi}\int_{0}^{2\pi} \frac{e^{it}}{(e^{it}-r\eit)^2}\, u(e^{it})\, dt.
\ee
Also, from the Cauchy integral formula, it is easy to see
$$0= \frac{1}{2\pi i}\int_{\IT} \frac{d\zeta}{(\zeta-z)^2} =  \frac{1}{2\pi}\int_{0}^{2\pi} \frac{e^{it}}{(e^{it}-r\eit)^2}\, dt,$$ so that \be\label{eq02}
0=\frac{1}{\pi} \int_{0}^{2\pi} \frac{e^{it}}{(e^{it}-r\eit)^2}\, u(\eit)\, dt.
\ee
Subtracting \eqref{eq02} from \eqref{eq01}, and taking absolute value, we find \be\label{eq03}|F'(r\eit)| \le \frac{1}{\pi} \int_{0}^{2\pi} \frac{\left|u(e^{i(\theta+t)})-u(\eit)\right|}{1-2r\cos t +r^2} \, dt.\ee Since $u(\eit)\in \Lambda_{\alpha}$, we have $$\left|u(e^{i(\theta+t)})-u(\eit)\right| \le A|t|^\alpha,$$ for some constant $A>0$.
Therefore, it follows from \eqref{eq03} that $$|F'(r\eit)| \le \frac{A}{\pi} \int_{0}^{2\pi} \frac{|t|^\alpha}{1-2r\cos t +r^2} \, dt =  \frac{2A}{\pi} \int_{0}^{\pi} \frac{t^\alpha}{1-2r\cos t +r^2} \, dt.$$ For $0\le t \le \pi$, we can estimate the denominator as $$1-2r\cos t +r^2=(1-r)^2+4r\sin^2 \frac{t}{2} \le (1-r)^2+\frac{4r}{\pi^2}\,t^2,$$ which implies $$|F'(r\eit)| \le  \frac{2A}{\pi} \int_{0}^{\pi} \frac{t^\alpha}{(1-r)^2+(4r/\pi^2)t^2} \, dt.$$ Now, we substitute $u=t/(1-r)$ to obtain \begin{align*}
	|F'(r\eit)| & \le  \frac{2A}{\pi} \frac{1}{(1-r)^{1-\alpha}}\int_{0}^{\pi/(1-r)} \frac{u^\alpha}{1+(4r/\pi^2)u^2} \, dt\\ & \le \frac{2A}{\pi} \frac{1}{(1-r)^{1-\alpha}}\int_{0}^{\infty} \frac{u^\alpha}{1+(4r/\pi^2)u^2} \, dt\\ & = O\left(\frac{1}{(1-r)^{1-\alpha}}\right),\end{align*} because the last integral converges for $\alpha < 1$. As in the proof of Theorem \ref{thm1}, we have $$|h'| \le \frac{1}{1-k}\, |F'|, \quad |g'| \le \frac{k}{1-k}\, |F'|,$$ and therefore, $$|h'(r\eit)|=O\left(\frac{1}{(1-r)^{1-\alpha}}\right)=|g'(r\eit)|.$$ Then, an appeal to Theorem \ref{smooth} implies $$h(\eit) \in \Lambda_{\alpha} \quad \text{and} \quad g(\eit) \in \Lambda_{\alpha}.$$ It follows that $f(\eit) \in \Lambda_{\alpha}$, and consequently, $v(\eit) \in \Lambda_{\alpha}$, as desired. This completes the proof.
\epf
The theorem is not true for $\alpha=1$, even if $f$ is analytic (i.e., $1$-quasiregular). The following example is well-known.

\begin{exm}
Let \( u \) be the harmonic function  in \( \mathbb{D} \) with boundary values
\[
u(e^{i\theta}) = |\theta| \quad \text{for } \theta \in [-\pi, \pi].
\]
Clearly, \( u(\eit) \) is Lipschitz. One can show, by the method of Hilbert transforms, that the boundary values of the conjugate function \( v \) behave like
\[
v(e^{i\theta}) \sim \theta \log|\theta| \quad \text{near } \theta = 0.
\]
It follows that 
\[
v'(e^{i\theta}) \sim \log|\theta|,
\]
which is unbounded as \( \theta \to 0 \). Thus, \( v(\eit) \) is not Lipschitz.
\end{exm}

\brem
The H\"older continuity of quasiregular mappings have been widely studied in the literature. Suppose $G \subset \IR^n$, $n\ge 2$, is a domain and $\mathbb{B}^n$ is the unit ball in $\IR^n$. It is known (see \cite[Theorem 1.11]{Rickman_quasiregular}, cf. \cite[Theorem 16.13]{HKVbook}) that every bounded $K$-quasiregular mapping $f:G\to \IR^n$ is $\delta$-H\"older continuous for some exponent $\delta \in (0,1]$ which depends on the \textit{inner dilatation} of $f$ (and therefore, on the constant $K$). Further, the exponent $\delta$ is best possible, as can be seen from the function $f: \mathbb{B}^n \to \mathbb{B}^n$, $f(x)=|x|^{\delta-1}x$ (here $\delta=K^{1/(1-n)}$).

It is important to clarify that Theorem \ref{thm2} presented herein diverges from this setting. We have shown that if $u(\eit)$ is $\alpha$-H\"older continuous, then so is $v(\eit)$, for any arbitrary $\alpha\in (0,1)$, i.e., the constant $K$ plays no role here. In other words, the primary importance of our result is in showing that the real and imaginary parts of a (planar) harmonic quasiregular mapping essentially behave like ``harmonic conjugates".
\erem

\subsection*{Acknowledgements} The research was partially supported by the Li~Ka~Shing Foundation STU-GTIIT Joint Research Grant (Grant no. 2024LKSFG06) and the Natural Science Foundation of Guangdong Province (Grant no. 2024A1515010467).


\bibliography{references}


\end{document}